\documentclass[12pt]{article}
\usepackage{color}

\usepackage{amsmath, amssymb, theorem, latexsym, color}
\usepackage{graphicx}

\newcommand{\beq}{\begin{equation}}
\newcommand{\eeq}{\end{equation}}

{}
\def\Ab{{\bf A}}
\def\Om {{\Omega}}
\def\la {{\lambda}}

\def \Ab{{\bf A}}
\def \Bb {{\bf B}}
 
\def\Xb{{\bf X}}

\def \Inte{{\rm Int\,}}

\newcommand {\pa}{\partial}

\newtheorem{theorem}{Theorem}[section]

\newtheorem{proposition}[theorem]{Proposition}

\newtheorem{definition}[theorem]{Definition}
\newtheorem{remark}[theorem]{Remark}

\begin{document}

\centerline{\bf On a magnetic characterization of
 spectral minimal partitions.
 }

\centerline{}

\centerline{}
 \centerline{ B. Helffer (Universit\'e Paris-Sud 11) }
 \centerline{and}
 \centerline{
 and T. Hoffmann-Ostenhof (University of Vienna)}

\centerline{} \centerline{}

\begin{abstract}
Given a bounded  open set $\Omega$ in $ \mathbb R^n$ (or  
 in a Riemannian  manifold)  and a partition of
$\Omega$ 
by $k$
 open sets $D_j$, we  consider the quantity
 $\max_j \lambda(D_j)$
 where $\lambda(D_j)$ is the ground state energy of
 the
 Dirichlet
 realization
 of the Laplacian in $D_j$. If we denote by
 $ \mathfrak
 L_k(\Omega)$ the infimum over all the $k$-partitions
 of  $ \max_j
 \lambda(D_j)$,  a minimal $k$-partition
 is then a partition which realizes the infimum. 
When
 $k=2$, we find the two  nodal domains of  a second
eigenfunction, but 
 the analysis of higher $k$'s is  non trivial and quite
 interesting. 
In this paper, we give the proof of
 one conjecture formulated in \cite{BH} and \cite{HeEg} about a magnetic characterization of the minimal partitions when $n=2$.\end{abstract}
{\bf{Keywords:}} minimal partitions, nodal sets, Aharonov-Bohm Hamiltonians,  Courant's nodal theorem.
\newline {\bf{AMS Subject
classification:}} 35B05.
\centerline{} \centerline{}
\section{Introduction}\label{Section1}
\subsection{Main definitions}

We consider mainly the Dirichlet  Laplacian in a bounded
domain $\Omega\subset \mathbb R^2$.
We would like to analyze the  relations between the nodal domains
 of the real-valued  eigenfunctions of this Laplacian and the partitions of
 $\Omega$ 
 by $ k$
 open sets $ D_i$ which are minimal in the sense that  the
 maximum over the $ 
D_i$'s of the ground state energy\footnote{The ground state energy is
  the smallest eigenvalue.}  of  the Dirichlet
 realization
 of the Laplacian $H(D_i)$  in $ D_i$ is minimal. In the case of a
Riemannian compact manifold, the natural extension is to consider 
 the Laplace Beltrami operator.
  We denote by $ \lambda_j(\Omega)$ the
increasing sequence of its eigenvalues and by $u_j$ some associated
 orthonormal basis of real-valued eigenfunctions. 
The ground state  $ u_1$ can be chosen to be
strictly positive in $ \Om$, but the other eigenfunctions 
$  u_k$ must have zerosets.
For any real-valued 
$ u\in C_0^0(\overline\Om)$, we define the zero set as 
\begin{equation}
N(u)=\overline{\{x\in \Om\:\big|\: u(x)=0\}}
\end{equation}
and call the components of $ \Om\setminus N(u)$ the nodal
domains of $ u$.
The number of 
nodal domains of $ u$ is called $ \mu(u)$. These
$\mu(u)$ nodal
domains define a $k$-partition of $ \Omega$, with $k=\mu(u)$. 

We  recall that the  Courant nodal
 theorem says that, for  $k\geq 1$,  and if $\lambda_k$ denotes  the
$k$-th eigenvalue  and  $ E(\lambda_k)$ the eigenspace of
$ H(\Omega)$ associated with $\lambda_k$, then,   for all real-valued $ u\in  E(\lambda_k)\setminus \{0\}\;,\;
\mu (u)\le k\;.
$

In dimension $1$ the Sturm-Liouville theory says that we have
always equality (for Dirichlet in a bounded interval) in the previous theorem (this is what we will call
later a Courant-sharp situation).  A  theorem due to Pleijel \cite{Pleijel:1956} in
1956  says
that this cannot be true when the dimension (here we consider the
$2D$-case) is larger than one.

We now introduce for $k\in \mathbb N$ ($k\geq 1$),
 the notion of $k$-partition. We
will call {\bf  $ k$-partition}  of $ \Omega$ a family 
$ \mathcal D=\{D_i\}_{i=1}^k$ of mutually disjoint sets in $\Omega$.
We  call it {\bf open} if the $D_i$ are open sets of
$ \Omega$,
{\bf connected}
 if the $ D_i$ are connected. We denote by $ \mathfrak O_k(\Omega)$ the set of open connected
partitions of $\Omega$. 
We now introduce the notion of spectral minimal partition sequence.
\begin{definition}\label{regOm}~\\
For any integer $ k\ge 1$,
 and for $ \mathcal D$ in $ \mathfrak O_k(\Omega)$, we
introduce 
\begin{equation}\label{LaD}
\Lambda(\mathcal D)=\max_{i}\la(D_i).
\end{equation} 
Then we define
\begin{equation}\label{frakL} 
\mathfrak L_{k}(\Omega)=\inf_{\mathcal D\in \mathfrak O_k}\:\Lambda(\mathcal D).
\end{equation}
and  call  $ \mathcal D\in \mathfrak O_k$ a minimal $k$-partition if 
  $ \mathfrak L_{k}=\Lambda(\mathcal D)$. 
\end{definition}

If $ k=2$, it is rather well known
 (see  \cite{HH:2005a} or
 \cite{CTV:2005}) that  $ \mathfrak L_2 =\lambda_2$ and
 that the associated minimal  $ 2$-partition is
 a {\bf nodal partition}, i.e. a partition whose elements
 are  the nodal domains of
   some eigenfunction corresponding to
 $ \lambda_2$. 

A partition $ \mathcal D=\{D_i\}_{i=1}^k$ of  $
\Omega$ in $ \mathfrak
O_k$ is called {\bf strong} if 
\begin{equation}\label{defstr}
\Inte(\overline{\cup_i D_i}) \setminus \pa \Om =\Om\;,
\end{equation}
where, for a set $A\subset  \mathbb R^2$,  $\Inte (A)$ means the interior of $A$.

Attached to a strong  partition, we  associate a closed
set in $ \overline{\Omega}$, which is called the {\bf boundary set}  of the partition~:
\begin{equation}\label{assclset} 
N(\mathcal D) = \overline{ \cup_i \left( \partial D_i \cap \Omega
  \right)}\;.
\end{equation}
$ N(\mathcal D)$ plays the role
 of the nodal set (in the case of a nodal partition).

This suggests the following definition: 
\begin{definition}\label{AMS}~\\
We call a partition  $\mathcal D$ regular if its associated
 boundary set  $ N(\mathcal D) $, has the following properties~:\\
(i)
Except for finitely many distinct $ x_i\in \Om\cap N$
 in the neighborhood of which $ N$ is the union of $\nu_i= \nu(x_i)$
smooth curves ($ \nu_i\geq 3$) with one end at $ x_i$,  
$ N$ is locally diffeomorphic to a regular 
curve.\\
(ii)
$ \pa\Om\cap N$ consists of a (possibly empty) finite set
of points $ z_i$. Moreover  
 $N$ is near $ z_i$ the union 
of $ \rho_i$ distinct smooth half-curves which hit
$ z_i$.\\
(iii) $ N$  has the {\bf equal angle
  meeting
 property}
\end{definition}
The $x_i$ are called the critical points and define the set
$X(N)$. Similarly
 we denote by $Y(N)$ the set of the boundary points $z_i$. 
By {\bf equal angle meeting property}, we mean that   the half curves meet with equal angle at each critical
 point of $ N$ and also at the boundary together with the
 tangent to the boundary.

 We say that $ D_i,D_j$ are {\bf  neighbors}
or $ D_i\sim D_j$,  if $
D_{ij}:=\Inte(\overline {D_i\cup D_j})\setminus \pa \Om$ is
connected. 
We associate with  
each $ \mathcal D$ a {\bf graph}
  $ G(\mathcal D)$ by
associating with  each $ D_i$ a vertex and to each 
pair $ D_i\sim D_j$ an edge. We will say that the graph is
{\bf bipartite} if it
can be colored by two colors (two neighbours having two different
colors). We recall that the graph associated
 with  a collection of nodal domains of an eigenfunction is always
 bipartite.
 
 \subsection{Motivation and outlook}
Before we state some results on spectral minimal partitions, discuss   their 
properties and finally formulate and prove the central result of the present paper, we give an informal outlook on our results.
The main result is a new characterization of minimal partitions via specific magnetic Hamiltonians, see Section \ref{Section4} for 
the necessary definitions and explanations of those operators. 
 
In \cite{HHOT} we have characterized via minimal partitions the case of equality in Courant's nodal theorem, 
see Theorem \ref{L=L} below. Roughly speaking, see Theorem \ref{partnod}, if a minimal partition could in principle stem 
from an eigenfunction it must be already be produced by the nodal domains of an eigenfunction and this can only happen if
there is equality in \eqref{Courant}. Pleijel's  result, \cite{Pleijel:1956}, implies, roughly speaking, 
that eigenfunctions associated to higher eigenvalues cannot lead to equality in \eqref{Courant}.  

In Section \ref{Section3} we  give a few pictures of non-nodal minimal partitions, or more precisely natural candidates, 
since it is notoriously hard to work out explicit examples for such partitions. A first glance shows that 
there are points where an odd number of nodal arcs meet. 

More than 10 years ago  together with Maria Hoffmann-Ostenhof 
and Mark  Owen we investigated some special magnetic Schr\"odinger operators, called Aharonov Bohm Hamiltonians, i.e.
 Hamiltonians with zero magnetic field but with singular
 magnetic vector potential and  with half integer circulation around holes in \cite{HHOO, HHOO1}, see Section \ref{Section4}.
 This investigation was  motivated by the at this time surprising result of Berger and Rubinstein, \cite{BeRu}, about the 
zeroset of a groundstate for such a problem with one hole. For more than one hole similar results were obtained 
on zerosets: each hole was hit by an odd number of nodal arcs.\footnote{ Similar results for 
punctured domains were later obtained in \cite{AFT}.} 

The findings in \cite{HHOO, HHOO1} motivated the conjecture in \cite{BH} and \cite{HeEg} and is reformulated in the present paper.  
The result says roughly that spectral minimal partitions are obtained by minimizing a certain eigenvalue 
of a Aharonov Bohm Hamiltonian with respect to the number and the position of poles if we assume that $\Om$ is simply connected.
See Theorem \ref{Theorem5.1} for the full result.

This new approach to spectral minimal partitions sheds new light on those spectral minimal partitions. While in in original
formulation, \cite{HHOT}, say for a fixed $\Om$ the $\mathfrak L_k(\Om)$ and the associated minimal partitions as defined by 
Definition \ref{regOm} require the calculation of $\Lambda(\mathcal D)$ for k-partitions, the new formulation can be considered as an, 
admittedly involved, eigenvalue minimization.

\paragraph{Acknowlegments}~\\ 
When writing this paper we benefitted from useful discussion with V. Bonnaillie-No\"el and S. Terracini.

 \section{Basic properties of  minimal partitions}\label{Section2}
 The following theorem has been proved  by Conti-Terracini-Verzini \cite{CTV0, CTV2,
  CTV:2005}
 and  Helffer--T.~Hoffmann-Ostenhof--Terracini \cite{HHOT}:
\begin{theorem}\label{thstrreg}~\\
For any $ k$, there exists a  minimal  regular $
k$-partition. Moreover 
any  minimal $ k$-partition has a   regular 
representative\footnote{Modulo sets of capacity $0$.}.
\end{theorem}
Other proofs of a  somewhat weaker version of this statement have been
given by Bucur-Buttazzo-Henrot \cite{BBH}, Caffarelli- F.H. Lin \cite{CL1}.

A natural question is whether a minimal partition of $ \Omega$ 
 is a  nodal partition, i.e. the family of
 nodal domains of  an eigenfunction of $ H(\Omega)$.
 We have  first the following converse theorem (\cite{HH:2005a}, \cite{HHOT}):
\begin{theorem}\label{partnod}~\\
If the graph of a minimal partition
 is bipartite,  then  this partition  is  nodal.
\end{theorem}
A natural question is now to determine how general   the previous 
situation is. Surprisingly this  only occurs in the so
 called Courant-sharp situation.  
  We say that $ u$ is {\bf  Courant-sharp} if
 $$  u\in  E(\lambda_k)\setminus \{0\} \quad \mbox{
  and}\quad 
 \mu(u)=k\;.
$$
For any integer $ k\ge 1$, we denote by $ L_k(\Omega)$
 the smallest eigenvalue of $H(\Omega)$, whose eigenspace 
contains an eigenfunction with $ k$ nodal domains.  We set 
$L_k(\Omega) =\infty$,  if there are no 
eigenfunction with $ k$ nodal domains.
  In general, one can 
 show that 
\begin{equation} 
\lambda_k(\Omega) \leq\mathfrak L_k(\Omega) \leq L_k(\Omega) \;.
\end{equation}
The last  result gives the full picture of the equality cases~:
\begin{theorem}\label{L=L}~\\
Suppose $ \Om\subset \mathbb R^2$ is regular.
If $\mathfrak L_k(\Omega)=L_k(\Omega)$ or $\mathfrak L_k(\Omega)=\lambda_k(\Omega)$
then 
\begin{equation}\label{Courant}
 \la_k(\Omega)=\mathfrak L_k(\Omega)=L_k(\Omega)\,.
 \end{equation}
 In addition, one can find in $ E(\la_k) $ a Courant-sharp 
 eigenfunction.
\end{theorem}
This answered a question posed in \cite{BHIM} (Section 7). 
\begin{remark}~\\
Very recently spectral partitions for discrete problems, namely
quantum graphs, have been investigated in   \cite{Bandetal}. 
 \end{remark}

\section{Examples of minimal $k$-partitions for special domains}\label{Section3}
Using Theorem  \ref{L=L}, 
 it is now easier to analyze the situation for the disk or for rectangles
 (at least in the irrational case), since we have  just to check for which eigenvalues
 one can find associated Courant-sharp eigenfunctions. 
 
 The possible  topological types of a minimal partition $\mathcal D$   rely  essentially 
 on Euler's formula and the fact that the $D_i$'s have to be  nice, that means
\begin{equation}\label{nice}
\Inte(\overline D_i)\cap \Omega =D_i\,.
\end{equation}
Figures \ref{fig.5part} and \ref{fig.disk} illustrate possible situations.
 
\begin{proposition}\label{Euler}~\\
Let $U $ be an open set in $\mathbb R^2$ with piecewise-$C^{1}$ boundary
 and let  $N$ a closed set such that $U \setminus N$ has $k$ 
 components and such that $N$ satisfies the properties of Definition
 \ref{AMS}.  Let $b_0$ be the number of components of 
$\pa U$ and $b_1$ be the number of components of $N\cup\pa U$. Denote by $\nu(x_i)$ and $\rho(z_i)$
the numbers of arcs associated with the $x_i\in X(N)$, respectively  $z_i\in Y(N)$. Then
\begin{equation}\label{Emu}
k =b_1-b_0+\sum_{x_i\in X(N)}(\frac{\nu(x_i)}{2}-1)+
\frac{1}{2}\sum_{z_i\in Y(N)}\rho(z_i)+1\,.
\end{equation}
\end{proposition} 
 This allows us  to analyze minimal partitions of a specific 
 topological type. If in addition the domain has some symmetries and we assume that a minimal partition keeps some of these symmetries,  then we find 
 natural  candidates for minimal partitions.
\paragraph{Minimal $3$-partitions}~\\
In the case of the disk (see \cite{HH:2006}), we have no proof that the minimal $ 3$-partition
 is the ``Mercedes star'' or $Y$-partition,  i.e. the partition created by three straight rays meeting
  at the center with equal angle.  But if   we assume that the minimal
 $ 3$-partition has a unique singular point at the center then one can show that  is indeed the $Y$-partition.This point of view is explored numerically 
by Bonnaillie-Helffer \cite{BH} (using some method equivalent to the Aharonov-Bohm approach and playing with the location of the critical point). There is also an interesting  theoretical analysis 
by Noris-Terracini \cite{NT}.\\
We have no example of minimal $3$-partitions with two critical points.
For the disk and the square the minimal $4$-partitions are nodal.

\paragraph{Minimal $5$-partitions}~\\
Using the covering approach, we were able (with V. Bonnaillie) in \cite{BH}  to produce numerically the
following
 candidate $\mathcal D_1$ for a minimal $ 5$-partition assuming  a specific topological
 type.
\begin{figure}[h!bt]
\begin{center}
\includegraphics[height=3cm]{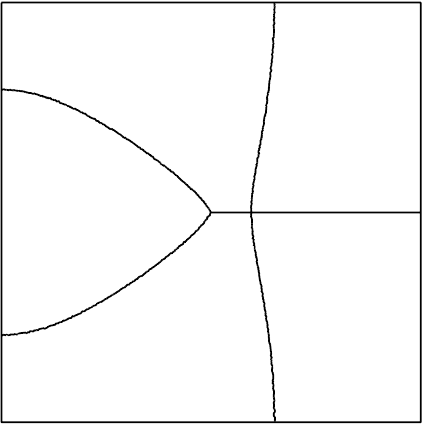}
\caption{Candidate $\mathcal D_1$ for the $5$-partition of the square.}
\end{center}
\end{figure}

It is interesting to compare with other possible topological types  of
  minimal $5$-partitions. They can be classified by using Euler's
  formula (see formula  \eqref{Emu}). Inspired by numerical
  computations in \cite{CyBaHo}, one looks for a 
 configuration which has the symmetries of the square and four
 critical points. We get two types of models that we can reduce
 to a Dirichlet-Neumann problem on a triangle corresponding to the
 eighth of the square. Moving the Neumann boundary on one side 
 like in \cite{BHV} leads us  to two candidates $\mathcal D_2$ and
  $\mathcal D_3$.  
One has a lower energy $\Lambda(\mathcal D)$ and
 one recovers  the pictures in \cite{CyBaHo}.
\begin{figure}[h!bt]
\begin{center}
\begin{tabular}{ccc}
$\Lambda(\mathcal D_1)=   111.910$
& $\Lambda(\mathcal D_2)=104.294$
& $\Lambda(\mathcal D_3) =131.666$\\
\includegraphics[height=3cm]{carrePart5}
& \includegraphics[height=3cm]{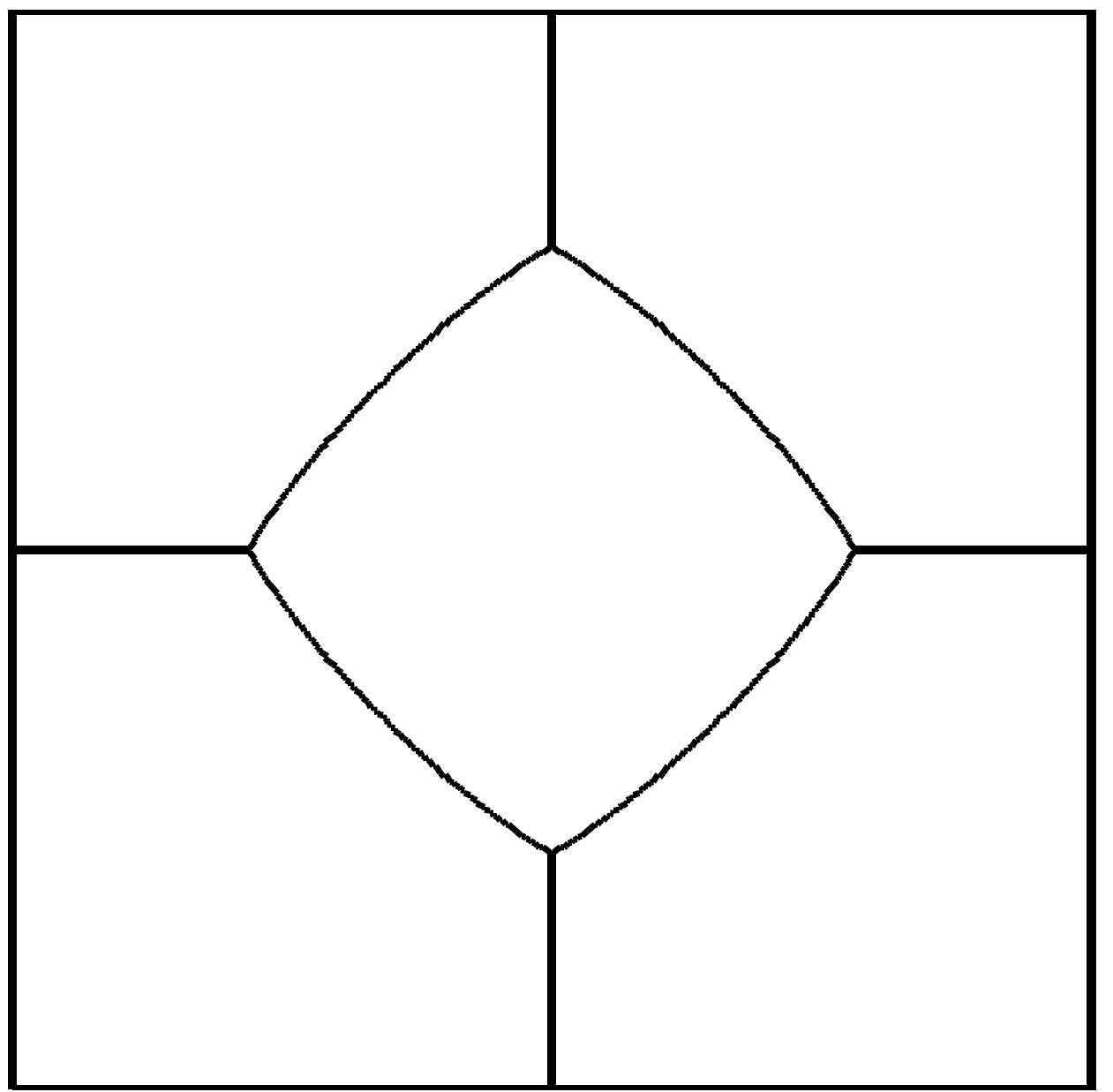}
& \includegraphics[height=3cm]{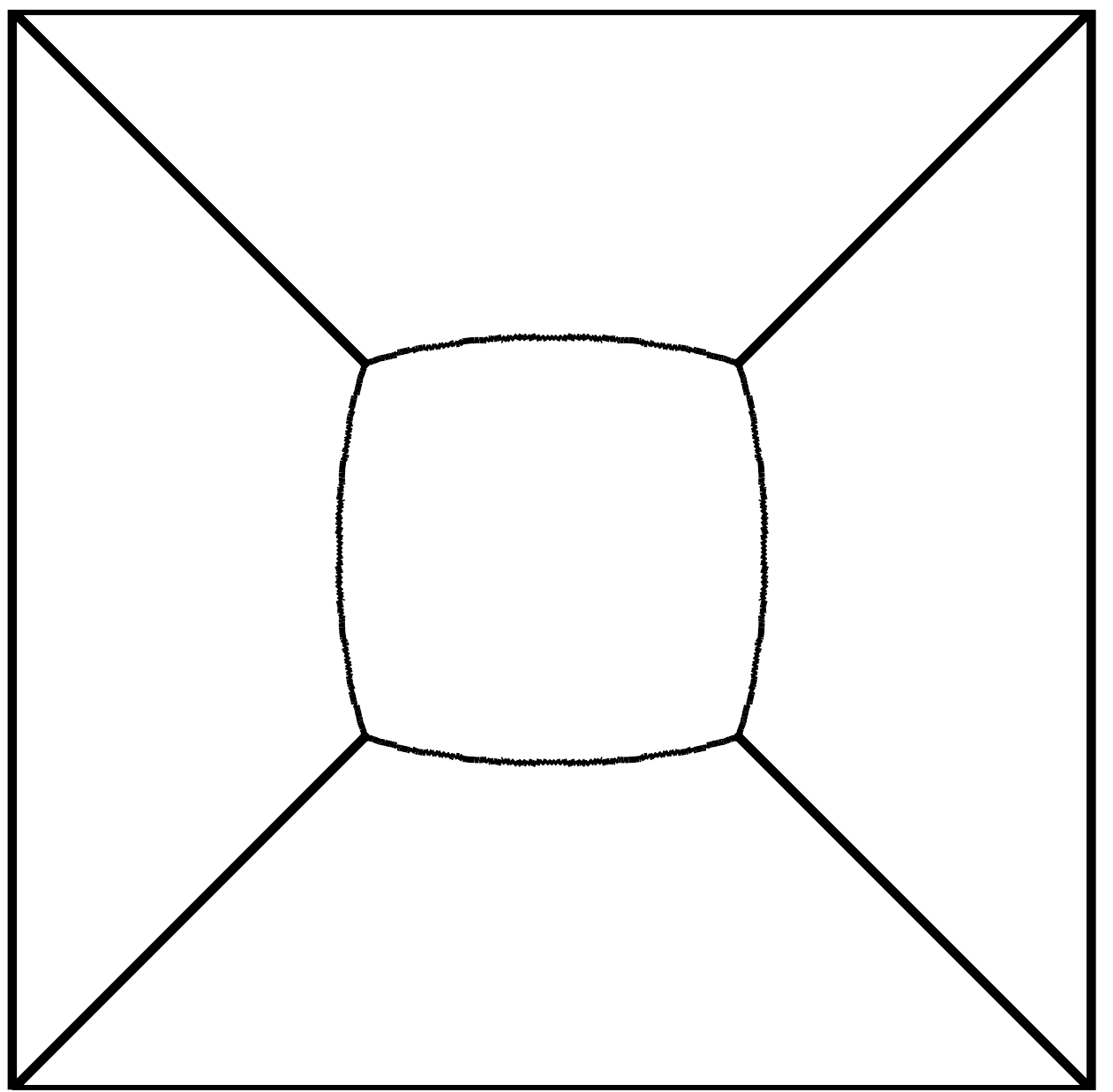}
\end{tabular}
\caption{Three  candidates for the $5$-partition of the square.\label{fig.5part}}
\end{center}
\end{figure}

Note that in the case of the disk a similar analysis leads to a
different answer. The partition of the disk by five half-rays with equal
angles has a lower energy than the
minimal $ 5$-partition with 
four singular points.

\begin{figure}[h!bt]
\begin{center}
\begin{tabular}{cc}
$104.367$
& $110.832$\\
\includegraphics[height=3cm]{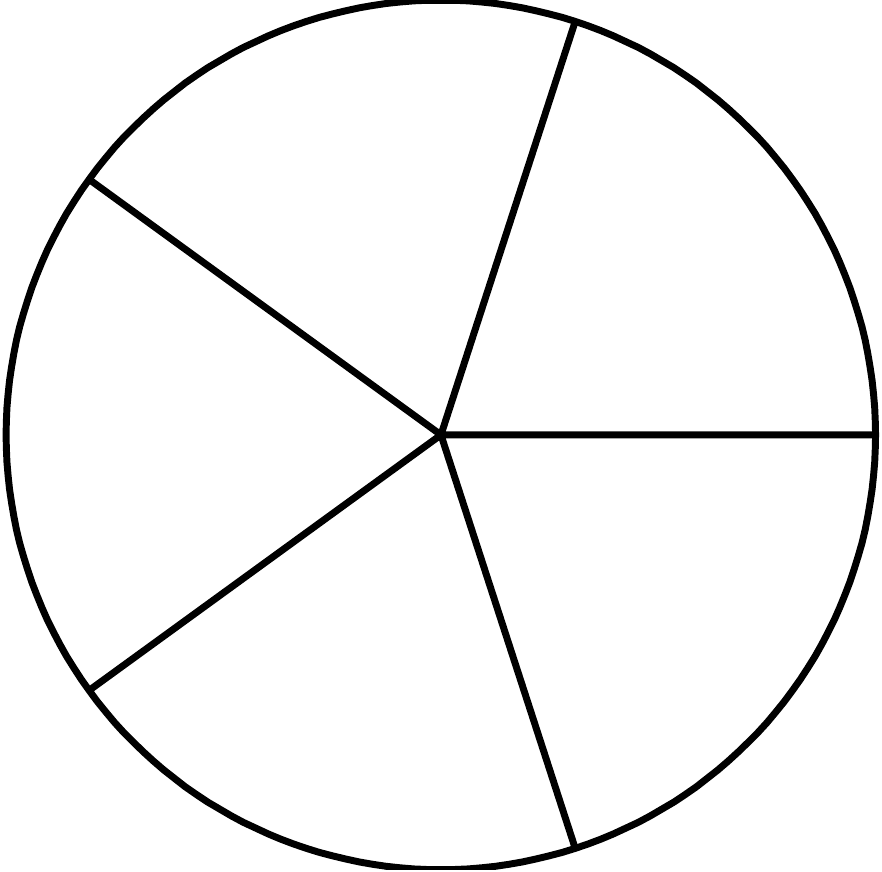}
& \includegraphics[height=3cm]{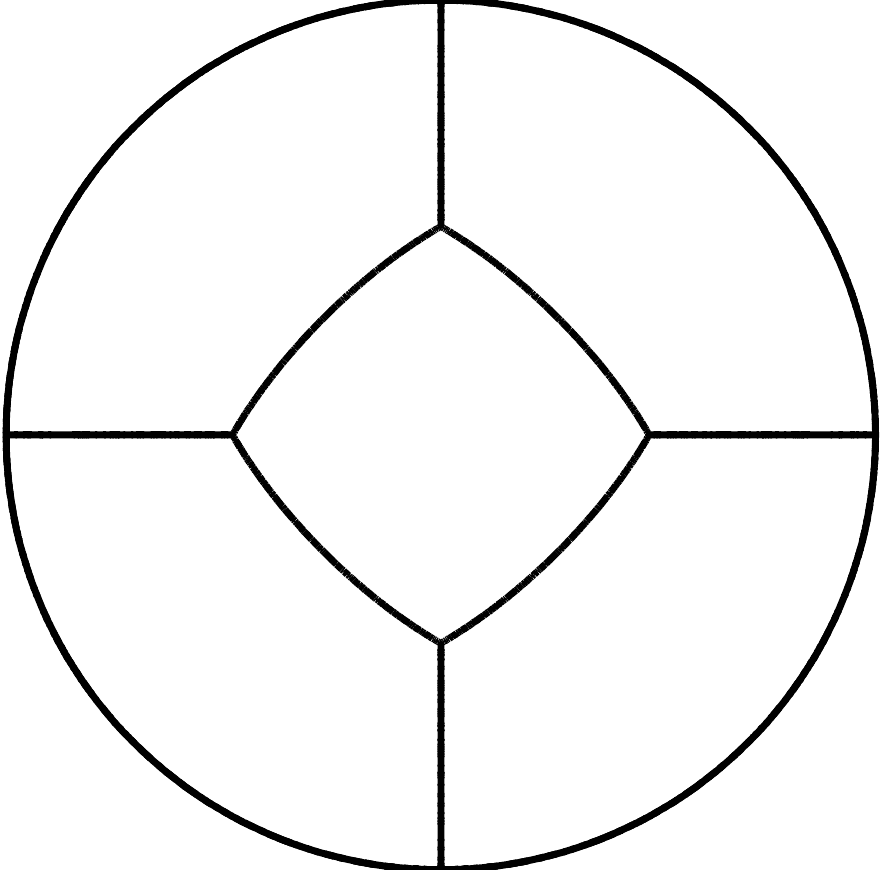}
\end{tabular}
\caption{Two candidates for the $5$-partition of the disk.\label{fig.disk}}
\end{center}
\end{figure}

\section{The Aharonov-Bohm approach}\label{Section4}

Let us recall some definitions and results about the Aharonov-Bohm 
Hamiltonian (for short ${\Ab\Bb}X$-Hamiltonian)  defined in an open set $\Omega$ which can be simply connected or not. These results were initially  motivated 
by the work of Berger-Rubinstein \cite{BeRu}, and further developed in  \cite{AFT, HHOO, HHOO1,BHHO,BH}.  
\paragraph{Simply connected case : one pole}~\\
We first consider the case when one pole, denoted by $X=(x_{0},y_{0})$, is chosen  in $\Omega$ and 
introduce the magnetic potential~:
\begin{equation}
{{\bf A}^X}(x,y) = (A_1^X(x,y),A_2^X(x,y))=\frac{\Phi}{2\pi}  \, \left( -\frac{y-y_{0}}{r^2}, \frac{x-x_{0}}{r^2}\right)\,.
\end{equation}
We know that in this case the magnetic field vanishes identically in
  $\dot\Omega_{X}\,$, where
  \begin{equation}
  \dot \Omega_{X}= \Omega \setminus \{X\}\,.
  \end{equation}
The ${\Ab\Bb}X$-Hamiltonian is defined  by considering the Friedrichs
extension starting from $C_0^\infty(\dot \Omega_{X})$
 and the associated differential operator is
\begin{equation}
-\Delta_{{\bf A}^X} := (D_x - A_1^X)^2 + (D_y-A_2^X)^2\,\mbox{ with }D_x =-i\pa_x\mbox{ and }D_y=-i\pa_y.
\end{equation}
We will consider in the sequel the very special case when the flux $\Phi$ created at $X=(x_0,y_0)$, which can be computed by considering the circulation of $\Ab^X$ along a simple closed path turning once anti-clockwise around $X$, satisfies:
\begin{equation}\label{condflux}
\frac{\Phi}{2\pi}=\frac 12 \,.
\end{equation}
Under assumption \eqref{condflux}, let $K_{X}$ be the anti-linear operator 
$$ K_{X} = e^{i \theta_{X}} \; \Gamma\,,$$
with $ (x-x_0)+ i (y-y_0) = \sqrt{|x-x_0|^2+|y-y_0|^2}\, e^{i\theta_{X}}\,$, 
where $\Gamma$ is the complex conjugation operator  
$$\Gamma u = \bar u\,
$$
and
\begin{equation}\label{deftheta}
 \nabla  \theta_X = 2 \Ab^X\,,
\end{equation}
which can also be rewritten in the form
$$
-\Ab^X = \Ab^X - \nabla  \theta_X\,.
$$
The flux condition \eqref{condflux} shows that one can find  a solution $\theta_X$ of \eqref{deftheta} (a priori multi-valued) such that   $e^{i \theta_{X}} $ is uni-valued and $C^\infty$.  Hence  $-\Delta_{{\bf A}^X}$ and  $-\Delta_{-{\bf A}^X}$ are intertwined by 
 the gauge transformation associated with $e^{i \theta_{X}} $. \\
Then we have
\begin{equation}\label{commuterel}
K_X \; \Delta_{{\bf A}^X}= \Delta_{{\bf A}^X}\; K_X\,.
\end{equation}
 We say that a function $u$ is  $K_{X}$-real, if it satisfies
$K_{X} u =u.$
Then the operator $-\Delta_{\Ab^X}$ is preserving the 
$K_{X}$- real functions. In the same way  one proves that the usual Dirichlet Laplacian admits an orthonormal  basis of real valued eigenfunctions or one restricts this Laplacian to the vector space over $\mathbb R$ of the real-valued $L^2$ functions,  one  can construct for  $-\Delta_{{\bf A}^X}$  a
basis of $K_{X}$-real eigenfunctions or, alternately, consider the 
restriction  of the ${\Ab\Bb}X$-Hamiltonian
to the vector space over $\mathbb R$ 
$$
L^2_{K_{X}}(\dot{\Omega}_{X})=\{u\in L^2(\dot{\Omega}_{X}) \;,\; K_{X}\,u =u\,\}\,.
$$
\paragraph{Non simply connected case}~\\
In this situation, magnetic potentials in $\Omega$ with zero magnetic field can be different from gradients if some fluxes around some holes are not in $(2\pi)\mathbb Z$. In this situation we will be interested in potentials where the created flux by some hole is $\pi$. This will be realized in this article by introducing a pole in the hole. Except that $\dot {\Omega}_X =\Omega$ (there are no singularity in $\Omega$) all what has been defined before goes through and this is actually the initial case treated in the pioneering work by \cite{BeRu}.\\

\paragraph{Poles and holes}~\\
We can extend our construction of an Aharonov-Bohm Hamiltonian
  in the case of a
 configuration with $\ell$ distinct points $X_1,\dots, X_\ell$ (putting a flux $\pi$ at each
 of these points). These points can be chosen in $\Omega$ or in the holes. They are distinct and each hole contains at most one $X_k$. We can just take as magnetic potential 
$$
\Ab^{\Xb} = \sum_{j=1}^\ell \Ab^{X_j}\,,
$$
where $\Xb=(X_1,\dots,X_\ell)$. 
Our Hamiltonian will be defined in
$
\dot{\Omega}_{\Xb} = \Omega \setminus \Xb\,.
$
 We can also construct (see \cite{HHOO,HHOO1}) the anti-linear
 operator $K_{\bf X}$,  where $\theta_X$ is replaced by a
 multivalued function $\phi_\Xb$ such that $\nabla \phi_\Xb = 2 \Ab^{\Xb}$ and $e^{i
 \phi_\Xb}$ is uni-valued and  $C^\infty$. We can then  consider the 
 real subspace of the $K_{\Xb}$-real
 functions in $L^2_{K_{\Xb}}(\dot{\Omega}_{\Xb})$ and our operator as an unbounded selfadjoint operator on $L^2_{K_{\Xb}}(\dot{\Omega}_{\Xb})$.

It was shown in \cite{HHOO,HHOO1} for the case with holes and in \cite{AFT} for the case with poles  that the nodal set of such a $K_X$-real eigenfunction has
the same structure as the nodal set of a real-valued  eigenfunction of the
Laplacian except that an odd number of half-lines meet at  each pole  and at the boundary of each hole containing some $X_k$.
 In the case of one hole, this fact was first observed by Berger-Rubinstein \cite{BeRu}  for a first
 eigenfunction (assuming that the first eigenvalue is simple). 
We denote
 by  $L_k(\dot{\Omega}_{\Xb})$ the lowest eigenvalue, 
 if it exists, such that there exists a $K_{\Xb}$-real eigenfunction with
 $k$ nodal domains and we set $L_k(\dot{\Omega}_{\Xb})=+\infty$ if there is no such eigenvalue.

\section{The magnetic
  characterization
 of a minimal partition}\label{Section5}
We now prove  the following conjecture presented (in the simply-connected case) in
 \cite{BH} and \cite{HeEg}. 
 \begin{theorem}\label{Theorem5.1}~\\
 Suppose $\Om$ is a bounded, not necessarily simply connected,  domain with $m$ disjoint closed holes $B_i$ ($i=1,\dots,m$) with non empty interiors. 
 Again we assume that $\pa\Om$ is piecewise $C^{1}$. Then 
\begin{equation} \label{magL}
 \mathfrak L_k(\Om)=\inf_{\ell\in \mathbb N}\:\inf_{X_1,\dots,X_\ell}L_k(\dot{\Om}_{\Xb})
\end{equation}
 where in the infimum   each $X_j= (x_j,y_j)$ is either in   $ \Inte(B_i)$ or in $\Om$. In each $B_i$ there is either one or no $X_i$. The $X_i \in \Om$
are distinct points.  
\end{theorem}

Let us first give the proof in the simply connected case. \\

{\bf Step 1} : $\inf_{\ell \in \mathbb N}\; \inf_{X_1,\dots,
  X_\ell} L_k (\dot{\Omega}_{\Xb})\,\leq \mathfrak L_k(\Omega)$\\
Considering a minimal $k$-partition $\mathcal D=(D_1,\dots,D_k)$, we know that it has a
regular representative and  we denote by $X^{odd}(\mathcal D):=(X_1,
\dots, X_\ell)$ the critical points of the boundary set of the partition for which 
an  odd number of half-curves meet. 

 For proving Step 1, we have indeed just to prove that, for this family of points $\Xb=X^{odd}(\mathcal D) $, 
 $\mathfrak L_k(\Omega)$ is an eigenvalue of the Aharonov-Bohm Hamiltonian associated with $\dot{\Omega}_{\Xb}$ 
  and to explicitly construct the corresponding eigenfunction with $k$ nodal domains described by the $D_i$'s.
 
For this, we recall that we have proven in
 \cite{HHOT} the existence of a family $(u_i)_{i=1,\dots,k}$ such that $u_i$ is a
 ground state of $\,H(D_i)$ and $u_i -u_j$ is a second eigenfunction of $H(D_{ij})$
 when $D_i\sim D_j$.  The claim  is that one can  find a sequence $\epsilon_i(x)$
 of $\mathbb S^1$-valued functions, where $\epsilon_i$ is a suitable\footnote{
   Note that by construction the $D_i$'s never contain any point of $\Xb$. Hence the ground state energy of the Hamiltonian $H(D_I)$ is the same 
    as the ground state energy of $H_{\Ab^\Xb} (D_i)$.}
 square root of  $e^{i\phi_{\Xb}}$ in $D_i$,
 such that
 $\sum_i \epsilon_i(x) u_i(x)$ is an eigenfunction of the $\Ab\Bb\Xb$-Hamiltonian associated with 
 the
 eigenvalue $\mathfrak L_k$.\\
 
 More explicitly, let us describe how we can construct $\epsilon_i(x)$. We start from some $i_0$ and define $\epsilon_{i_0} (x) = e^{\frac {i}{2}\phi_{\Xb}}$. According to the footnote
  $\epsilon_{i_0}(x)$ is a well defined $C^\infty$ function. Let $D_i$ a nearest neighbor of $D_{i_0}$ then we define $\epsilon_i(x) = - e^{\frac {i}{2} \phi_{\Xb}}$. Then we can 
   extend the definition by considering the neighbors of the neighbors. Now we have to check that the construction is consistent. The problem can be reduced to the following question. Consider a closed simple path $\gamma$ in $\dot{\Omega}_X$ transversal to $\mathcal N(\mathcal D)$ (and avoiding the critical points). Take some origin $x_0$ on $\gamma \cap D_{i_1}$. We start
    from $\epsilon(x) = e^{\frac i2 \phi_{\Xb}(x)}$ in $D_{i_1}$ and, choosing the positive orientation,  multiply by $-1$ each time that we cross an arc of $\mathcal N(\mathcal D)$. It is then a consequence of Euler's formula
     that the number of crossings along $\gamma$ is odd if and only if there is an odd number of points of $\Xb$ inside $\gamma$ (apply Euler's formula \eqref{Emu} with $U$ being the open set delimited by $\gamma$).  It is then clear that $\epsilon (x)$ is well defined along $\gamma$.\\
 
{\bf Step 2}: $\inf_{\ell \in \mathbb N}\; \inf_{X_1,\dots,
  X_\ell} L_k (\dot{\Omega}_{\Xb})\,\geq \mathfrak L_k(\Omega)$\\
Conversely, given $\ell$ distinct points $X_i$ in $\Omega$,  any family of nodal domains of a $K_{\Xb}$-real eigenfunction of  the Aharonov-Bohm
operator on $\dot{\Omega}_{\Xb}$ corresponding to $L_k$
 gives a $k$-partition.  Using the results of \cite{HHOO} and \cite{AFT}, we immediately see that the $X_i$'s corresponding to the "odd" singular points of the partitions. 
  In each of these nodal domains $D_i$,  $L_k$ is an eigenvalue of the Dirichlet realization of the Schr\"odinger operator with magnetic potential $\Ab^\Xb$, which is by the diamagnetic inequality higher as the ground state energy of the Dirichlet Laplacian in $D_i$ without magnetic field. Hence the energy $\Lambda_k (\mathcal D)$ of this partition is indeed less than $L_k(\dot{\Omega}_{\Xb}) $.\\
  
 {\bf Step 3}: Proof in the non simply connected case~\\
  The main change is in step 1. In the non simply connected case, the set $\Xb$ consists of the singular points of the boundary set inside $\Omega$ where an odd number of half-lines arrive  together with  those points  in the holes  whose  boundary is hit  by an odd number of half-curves.\\

 \paragraph{Examples}~\\
  Let us present a few examples  illustrating  the theorem in the case of a  simply connected domain. 
 When
$k=2$, there is no need to consider punctured $\Omega$'s. The infimum
 is obtained for $\ell =0$.  When $k=3$,
it is possible to show (see Remark \ref{Rem5.3} below) that it is enough to minimize over $\ell =0$,
$\ell =1$ and $\ell =2$. In the case of the disk and the square, it is 
 proven that the infimum cannot be for $\ell =0$ and we
conjecture that the infimum is for $\ell =1$ and attained for the punctured
domain  at the center. For $k=5$,  it seems
that the infimum is for $\ell =4$ in the case of the square (See Figure \ref{fig.5part})  and for $\ell =1$ in the case of the disk (see Figure \ref{fig.disk}).\\
 \begin{remark}\label{Rem5.2}~\\
  If $\mathcal D$ is a regular representative of a minimal $k$-partition and if $\dot{\Omega}_\Xb $ is constructed like in Step 1  of the proof of the previous theorem, then 
 $\mathfrak
 L_k(\Omega) =\lambda_k(\dot{\Omega}_\Xb)$ (Courant sharp
 situation). Coming back indeed to this step, one can follow the proof of Theorem 1.13 (Section 6) in \cite{HHOT}. 
 \end{remark}
 \begin{remark}\label{Rem5.3}~\\
Euler's formula \eqref{Emu}, implies that for a minimal $k$-partition
 $\mathcal D$ of a simply connected domain $\Omega$
 the cardinality of $X^{odd}(\mathcal D)$ satisfies
\begin{equation}
\# X^{odd}(\mathcal D ) \leq 2k -3\,. 
\end{equation}
Note that if $b_1=b_0$, we necessarily have a singular point in the boundary. The argument   depends only on Euler's formula.  If we implement the additional  property
 that the open sets $D_i$'s of a minimal partition are nice (see \eqref{nice}), we can exclude the case when there is only one point on the boundary. We emphasize that this was not a priori excluded from the results of \cite{HHOO,AFT}.  Hence, we obtain 
$$
b_1-b_0 + \frac 12 \sum \rho(y_i) \geq 1\,, 
$$
which implies the inequality 
\begin{equation}
\# X^{odd}(\mathcal D ) \leq 2k -4 \,. 
\end{equation}
This estimate seems optimal for a general geometry although all the known candidates for minimal partitions
 for $k=3$ and $5$ have  a lower cardinality of odd critical points. 
\end{remark}
\begin{remark}\label{Rem5.4}~\\
The argument around \eqref{nice} shows that a
nodal set of a $K_{\mathbf X}$-real eigenfunction that corresponds to a minimal partition
cannot have a critical point  that is met only by one nodal arc. Actually that can happen
for ground states of Aharonov-Bohm  Hamiltonians, see \cite{HHOO} which of course do not
correspond to minimal partitions.

\end{remark}

 \begin{remark}~\\
It would be interesting to look at the case of the sphere (already considered in \cite{HHOT1}) and the first problem in this case is to define the suitable magnetic Laplacian. 
We refer to \cite{WY} for one of the first papers on this question. More specifically, we would like to construct in our case an Aharomov-Bohm Hamiltonian.  Note for example that we can not have such an operator with one pole and a  flux $\pi$ around this pole. Fortunately there are no minimal $k$-partition  whose boundary set consists of one "odd" critical point  
 on the sphere, as can be seen by
  Euler's formula for the sphere (see in \cite{HHOT1}, Remark 4.2).  We indeed know that the cardinality of "odd" critical  points  is even.  This is actually 
  a standard result from graph theory that  the number of vertices with odd degree is even.
(See for example Corollary 1.2 in \cite{BM}).\\
This suggests that instead of putting the flux $\pi$ around each pole, we take
   alternately $\pi$ and $-\pi$ for the fluxes in order to get a total flux equal to $0$. In other words, we should probably describe $X^{odd} (\mathcal D)$ as a union of dipoles.
 \end{remark}

\scshape 
B. Helffer: D\'epartement de Math\'ematiques, Bat. 425,
Universit\'e Paris-Sud, 91 405 Orsay Cedex, France.

email: Bernard.Helffer@math.u-psud.fr\\

T. Hoffmann-Ostenhof: Department of Theoretical Chemistry, 1090 Wien,
W\"ahringerstrasse  17, Austria

email: thoffmann@tbi.univie.ac.at

\end{document}